\documentclass[a4paper,11pt,twoside]{amsart}

\usepackage[latin1]{inputenc}
\usepackage[T1]{fontenc}
\usepackage{amsmath}
\usepackage{amsthm}
\usepackage{amssymb}

\newtheorem{thm}{Theorem}[section]
\newtheorem{prop}[thm]{Proposition}
\newtheorem{lem}[thm]{Lemma}
\newtheorem{cor}[thm]{Corollary}

\theoremstyle{definition}

\newtheorem{remark}{Remark}[section]

\begin{document}

\title[FM-partners
of K3 surfaces with Picard numbers 1 and 2]{Some remarks about the
FM-partners of K3 surfaces with Picard numbers 1 and 2}

\author{Paolo Stellari}

\address{Dipartimento di Matematica ``F.
Enriques'', Universit{\`a} degli Studi di Milano, Via Cesare
Saldini 50, 20133 Milano, Italy}

\email{stellari@mat.unimi.it}



\begin{abstract} In this paper we prove some results about K3 surfaces with
Picard number 1 and 2. In particular, we give a new simple proof
of a theorem due to Oguiso which shows that, given an integer $N$,
there is a K3 surface with Picard number 2 and at least $N$
non-isomorphic FM-partners. We describe also the Mukai vectors of
the moduli spaces associated to the Fourier-Mukai partners of K3
surfaces with Picard number 1.\end{abstract}

\maketitle

\section{Introduction}\label{sec:intro}

In some recent papers Hosono, Lian, Oguiso and Yau (see
\cite{HLOY3} and \cite{Og}) gave a formula that counts the number
of non-isomorphic Fourier-Mukai partners of a K3 surface. In this
paper we are interested in the case of K3 surfaces with Picard
number 1 and 2.

In the second paragraph, we recall the formula for the number of
the isomorphism classes of Fourier-Mukai partners of a given K3
surface (given in \cite{HLOY2}), which allows to count the
isomorphism classes of Fourier-Mukai partners of a K3 surface with
Picard number 1 (this is also given in \cite{Og}). As a first
result, we will describe the Mukai vectors of the moduli spaces
associated to the Fourier-Mukai partners of such K3
surfaces\footnote{This result was independently proved by Hosono,
Lian, Oguiso and Yau (Theorem 2.1 in \cite{HLOY3}).}. This gives
some information about the geometry of the Fourier-Mukai partners
of the given K3 surface.

In the third paragraph we prove that, given $N$ and $d$ positive
integers, there is an elliptic K3 surface with a polarization of
degree $d$ and with at least $N$ non-isomorphic elliptic
Fourier-Mukai partners (Theorem \ref{thm:main2}). The most
interesting consequence of this result is a new simple proof of
Theorem 1.7 in \cite{Og} (Corollary \ref{cor:oguiso} and Remark
3.5).

We start with recalling some essential facts about lattices and K3
surfaces.


\subsection{Lattices and discriminant groups}

A {\it lattice} $L:=(L,b)$ is a free abelian group of finite rank
with a non-degenerate symmetric bilinear form $b:L \times L
\rightarrow\mathbb{Z}$. Two lattices $(L_1,b_1)$ and $(L_2,b_2)$
are {\it isometric} if there is an isomorphism of abelian groups
$f: L_1 \rightarrow L_2$ such that $b_1(x,y)=b_2(f(x),f(y))$. We
write $\mathrm{O}(L)$ for the group of all autoisometries of the
lattice $L$. A lattice $(L,b)$ is {\it even} if, for all $x\in L$,
$x^2:=b(x,x)\in 2\mathbb{Z}$, it is {\it odd} if there is $x\in L$
such that $b(x,x)\not \in 2\mathbb{Z}$. Given an integral basis
for $L$, we can associate to the bilinear form a symmetric matrix
$S_L$ of dimension rk$L$, uniquely determined up to the action of
$\mathrm{GL}(\mathrm{rk}L,\mathbb{Z})$. The integer
det$L$:=det$S_L$ is called {\it discriminant} and it is an
invariant of the lattice. A lattice is {\it unimodular} if
det$L$=$\pm1$. Given $(L,b)$ and $k\in\mathbb{Z}$, $L(k)$ is the
lattice $(L,kb)$.

Given a sublattice $V$ of $L$ with $V \hookrightarrow L$, the
embedding is {\it primitive} if $L/V$ is free. In particular, a
sublattice is primitive if its embedding is primitive. Two
primitive embeddings $V \hookrightarrow L$ and $V \hookrightarrow
L'$ are {\it isomorphic} if there is an isometry between $L$ and
$L'$ which induces the identity on $V$. For a sublattice $V$ of
$L$ we define the {\it orthogonal} lattice $V^\perp := \{x\in
L:b(x,y)=0, \forall y \in V\}$. Given two lattices $(L_1,b_1)$ and
$(L_2,b_2)$, their {\it orthogonal direct sum} is the lattice
($L,b$), where $L=L_1 \oplus L_2$ and $b(x_1+y_1,x_2+y_2)=
b_1(x_1,x_2)+b_2(y_1,y_2)$, for $x_1,x_2 \in L_1$ and $y_1,y_2 \in
L_2$.

The {\it dual lattice} of a lattice $(L,b)$ is
$L^\vee:=\mathrm{Hom}(L,\mathbb{Z})\cong \{x\in L \otimes
_\mathbb{Z}\mathbb{Q}: b(x,y)\in\mathbb{Z}, \forall y\in L\}$.
Given the natural inclusion $L\hookrightarrow L^\vee$, $x\mapsto
b(-,x)$, we define the {\it discriminant group} $A_L$:=$L^\vee/L$.
The order of $A_L$ is $\mid$det$L\mid$ (see \cite{BPV}, Lemma 2.1,
page 12). Moreover, $b$ induces a symmetric bilinear form $b_L:A_L
\times A_L \rightarrow \mathbb{Q}/\mathbb{Z}$ and a corresponding
quadratic form $q_L:A_L \rightarrow \mathbb{Q}/\mathbb{Z}$ such
that, when $L$ is even, $q_L(\overline{x})=q(x)$ modulo
$2\mathbb{Z}$, where $\overline{x}$ is the image of $x\in L^\vee$
in $A_L$. The elements of the triple $(t_{(+)},t_{(-)},q_L)$,
where $t_{(\pm)}$ is the multiplicity of positive/negative
eigenvalues of the quadratic form on $L\otimes\mathbb{R}$, are
invariants of the lattice $L$.

If $L$ is unimodular, $L^\vee\cong \{b(-,x):x\in L\}$. If $V$ is a
primitive sublattice of a unimodular lattice $L$ such that
$b\mid_V$ is non-degenerate, then there is a natural isometry of
groups $\gamma:V^\vee$/$V\rightarrow (V^\perp)^\vee/V^\perp$.


\subsection{K3 surfaces and $M$-polarizations}

A {\it K3 surface} is a 2-dimensional complex projective smooth
variety with trivial canonical bundle and first Betti number
$b_1=0$. From now on, $X$ will be a K3 surface. The group
$H^2(X,\mathbb{Z})$ with the cup product is an even unimodular
lattice and it is isomorphic to the lattice $\Lambda:=U^3\oplus
E_8(-1)^2$ (for the meaning of $U$ and $E_8$ see \cite{BPV} page
14). The lattice $\Lambda$ is called {\it K3 lattice} and it is
unimodular and even.

Given the lattice $H^2(X,\mathbb{Z})$, the {\it N\'eron-Severi
group} $\mathrm{NS}(X)$ is a primitive sublattice.
$T_X:=\mathrm{NS}(X)^\perp$ is the {\it transcendental lattice}.
The rank of the N\'eron-Severi group $\rho(X):=$rk$\mathrm{NS}(X)$
is called the {\it Picard number}, and the signature of the
N\'eron-Severi group is $(1,\rho -1)$, while the one of the
transcendental lattice is $(2,20-\rho)$. If $X$ and $Y$ are two K3
surfaces, $f:T_X\rightarrow T_Y$ is an {\it Hodge isometry} if it
is an isometry of lattices and the complexification of $f$ is such
that $f_\mathbb{C}(\mathbb{C}\omega_X)=\mathbb{C}\omega_Y$, where
$H^{2,0}(X)=\mathbb{C}\omega_X$ and
$H^{2,0}(Y)=\mathbb{C}\omega_Y$. We write
$(T_X,\mathbb{C}\omega_X)\cong (T_Y,\mathbb{C}\omega_Y)$ to say
that there is an Hodge isometry between the two transcendental
lattices.

A {\it marking} for a K3 surface $X$ is an isometry $\varphi:
H^2(X,\mathbb{Z})\rightarrow \Lambda$. We write $(X,\varphi)$ for
a K3 surface $X$ with a marking $\varphi$. Given
$\Lambda_\mathbb{C}:=\Lambda\otimes\mathbb{C}$ and given
$\omega\in \Lambda_\mathbb{C}$ we denote by
$[\omega]\in\mathbb{P}(\Lambda_\mathbb{C})$ the corresponding line
and we define the set
$\Omega:=\{[\omega]\in\mathbb{P}(\Lambda_\mathbb{C}):
\omega\cdot\omega=0, \omega\cdot \overline {\omega}>0\}$. The
image in $\mathbb{P}(\Lambda_\mathbb{C})$ of the line spanned by
$\varphi_\mathbb{C}(\omega_X)$ belongs to $\Omega$ and is called
{\it period point} (or period) of the marked surface
$(X,\varphi)$. From now on, the period point of a marked K3
surface $(X,\varphi)$ will be indicated either by
$\mathbb{C}\varphi_\mathbb{C}(\omega_X)$ or by
$[\varphi_\mathbb{C}(\omega_X)]$.

Given two K3 surfaces $X$ and $Y$, we say that they are {\it
Fourier-Mukai-partners} (or FM-partners) if there is an
equivalence between the bounded derived categories of coherent
sheaves $\mathrm{D}^b_{coh}(X)$ and $\mathrm{D}^b_{coh}(Y)$. By
results due to Mukai and Orlov, this is equivalent to say that
there is an Hodge isometry $(T_X,\mathbb{C}\omega_X)\rightarrow
(T_Y,\mathbb{C}\omega_Y)$. We define $FM(X)$ to be the set of the
isomorphism classes of the FM-partners of $X$.

Let $M$ be a primitive sublattice of $\Lambda$ with signature
$(1,t)$. A K3 surface $X$ with a marking
$\varphi:H^2(X,\mathbb{Z})\rightarrow\Lambda$ is a {\it marked
M-polarized K3 surface} if $\varphi^{-1}(M)\subseteq
\mathrm{NS}(X)$. A K3 surface is {\it $M$-polarizable} if there is
a marking $\varphi$ such that $(X,\varphi)$ is a marked
$M$-polarized K3 surface. Two marked and $M$-polarized surfaces
$(X,\varphi)$ and $(X',\varphi')$ are isomorphic if there is an
isomorphism $\psi:X\rightarrow X'$ such that
$\varphi'=\varphi\circ\psi^*$. Form now on, we will consider the
case of lattices $M:=\langle h\rangle$, with $h^2=2d$ and $d>0$.
The pair $(X,h)$, where $X$ is a K3 surface and $h\in
\mathrm{NS}(X)$, with $h^2=2d$, means a K3 surface with a
polarization of degree $2d$.


\section{FM-partners of a K3 surface with
$\rho=1$ and associated Mukai vectors}\label{sec:rk1}


In this section we want to describe the Mukai vectors of the
moduli spaces associated to the $M$-polarized FM-partners of a K3
surface $X$ with Picard number 1. By Orlov's results (\cite{Or}),
$q=|FM(X)|$ is the same as the number of non-isomorphic compact
2-dimensional fine moduli spaces of stable sheaves on $X$.
Obviously, on a K3 surface with Picard number 1 and
$\mathrm{NS}(X)=\langle h\rangle$ there is only one $\langle
h\rangle$-polarization of degree $h^2=2d$. So the concept of
FM-partner and the concept of $M$-polarized FM-partner coincide.
If $M=\langle h\rangle$ we are sure, by Orlov, that if we find $q$
non-isomorphic moduli spaces, then these are representatives of
all the isomorphism classes of $M$-polarized FM-partners of $X$.

We recall briefly the counting formula for the isomorphism classes
of FM-partners of a given K3 surface. Given a lattice $S$, the
{\it genus} of $S$ is the set ${\mathcal G}(S)$ of all the
isometry classes of lattices $S'$ such that $A_S\cong A_{S'}$ and
the signature of $S'$ is equal to the one of $S$.

Let $T_X$ be the transcendental lattice of an abelian surface or
of a K3 surface $X$ with period $\mathbb{C}\omega_X$. We can
define the group
$$
G:=O_{Hodge}(T_X,\mathbb{C}\omega_X)=\{g\in
\mathrm{O}(T_X):g(\mathbb{C}\omega_X)=\mathbb{C}\omega_X\}.
$$
We know (see \cite{Ca} Theorem 1.1, page 128), that the genus of a
lattice, with fixed rank and discriminant, is finite. The map
$\mathrm{O}(S)\rightarrow \mathrm{O}(A_S)$ defines an action of
$\mathrm{O}(S)$ on $\mathrm{O}(A_S)$. On the other hand, taken
$g\in G$, and given a marking $\varphi$ for $X$, $\varphi\circ
g\circ \varphi^{-1}$ induces an isometry on the lattice
$T:=\varphi(T_X)$, thus $\varphi$ defines a homomorphism
$G\hookrightarrow \mathrm{O}(T)$. The composition of this map and
the map $\mathrm{O}(T)\rightarrow \mathrm{O}(A_T)$ gives an action
of $G$ on $\mathrm{O}(A_T)\cong \mathrm{O}(A_S)$.

\begin{thm}\label{thm:form} {\bf [4, Theorem 2.3].} Let $X$ be a K3 surface and let
${\mathcal G}(\mathrm{NS}(X))={\mathcal G}(S)=\{S_1,\cdots,S_m\}$.
Then
$$
|FM(X)|=\sum^m_{j=1}|\mathrm{O}(S_j)\backslash
\mathrm{O}(A_{S_j})/G|,
$$
where the actions of the groups $G$ and $\mathrm{O}(S_j)$ are
defined as before.\end{thm}

The following corollary (which is Theorem 1.10 in \cite{Og})
determines the number $q$ of FM-partners of a surface with Picard
number 1.

\begin{cor}\label{cor:num} Let $X$ be a K3 surface with $\rho(X)=1$ and such that
$\mathrm{NS}(X)=\langle h\rangle$, with $h^2=2d$.

{\rm (i)} The group $\mathrm{O}(A_S)$ is trivial if $d=1$ while, if
$d>1$, $\mathrm{O}(A_S)\cong (\mathbb{Z}/2\mathbb{Z})^{p(d)}$, where
$p(d)$ is the number of distinct primes $q$ such that $q|d$. In
particular, if $d\geq 2$, then $|\mathrm{O}(A_S)|=2^{p(d)}$.

{\rm (ii)} For all markings $\varphi$ of $X$, the image of
$H_{X,\varphi}:=\{\varphi\circ g\circ \varphi^{-1}:g\in G\}\subseteq
\mathrm{O}(T)$ in $\mathrm{O}(A_T)$ by the map
$\mathrm{O}(T)\rightarrow \mathrm{O}(A_T)$ is $\{\pm
\overline{id}\}$.

In particular, $|FM(X)|=2^{p(d)-1}$, where now we set
$p(1)=1$.\end{cor}

Assertion (i) is known and it can also be found in \cite{Sc}
(Lemma 3.6.1).

Using the notation of \cite{Mu2}, we put
$H^*(X,\mathbb{Z}):=H^0(X,\mathbb{Z})\oplus H^2(X,\mathbb{Z})
\oplus H^4(X,\mathbb{Z})$. Given $\alpha:=(\alpha _1,\alpha
_2,\alpha _3)$ and $\beta:=(\beta _1,\beta _2,\beta _3)$ in
$H^*(X,\mathbb{Z})$, using the cup product we define the bilinear
form
$$
\alpha \cdot \beta:=-\alpha _1 \cup \beta _3+\alpha _2 \cup \beta
_2-\alpha _3 \cup \beta _1.
$$
From now on, depending on the context, $\alpha \cdot \beta$ will
mean the bilinear form defined above or the cup product on
$H^2(X,\mathbb{Z})$.

We give to $H^*(X,\mathbb{Z})$ an Hodge structure considering
\begin{eqnarray}
H^{*}(X,\mathbb{C})^{2,0} & := & H^{2,0}(X),\nonumber \\
H^{*}(X,\mathbb{C})^{0,2} & := & H^{0,2}(X),\nonumber \\
H^{*}(X,\mathbb{C})^{1,1} & := & H^0(X,\mathbb{C})\oplus
H^{1,1}(X)\oplus H^4(X,\mathbb{C}).\nonumber
\end{eqnarray}
$\tilde{H}(X,\mathbb{Z})$ is the group $H^*(X,\mathbb{Z})$ with
the bilinear form and the Hodge structure defined before.

For $v=(r,h,s)\in \tilde{H}(X,\mathbb{Z})$ with $r\in
H^0(X,\mathbb{Z})\cong\mathbb{Z}$, $s\in
H^4(X,\mathbb{Z})\cong\mathbb{Z}$ and $h\in H^2(X,\mathbb{Z})$,
$M(v)$ is the moduli space of stable sheaves $E$ on $X$ such that
rk$E=r$, $c_1(E)=h$ and $s=c_1(E)^2/2-c_2(E)+r$. If the stability
is defined with respect to $A\in H^2(X,\mathbb{Z})$ we write
$M_A(v)$. The vector $v$ is {\it isotropic} if $v\cdot v=0$. The
vector $v\in \tilde{H}(X,\mathbb{Z})$ is primitive if
$\tilde{H}(X,\mathbb{Z})/\mathbb{Z}v$ is free.

As we have observed, the results of Orlov in \cite{Or} imply that
each FM-partner of $X$ is isomorphic to an $M_h(v)$. We determine
a set of Mukai vectors which corresponds bijectively to the
isomorphism classes of the FM-partners of $X$ in $FM(X)$. First of
all, we recall the following theorem due to Mukai (\cite{Mu2}).

\begin{thm}\label{thm:spmod} {\bf [10, Theorem 1.5 3].} If $X$ is a K3 surface,
$v=(r,h,s)$ is an isotropic vector in
$\tilde{H}^{1,1}(X,\mathbb{Z})=H^{*}(X,\mathbb{C})^{1,1}\cap
H^*(X,\mathbb{Z})$ and $M_A(v)$ is non-empty and compact, then
there is an isometry $\varphi:v^\perp /\mathbb{Z}v\rightarrow
H^2(M_A(v),\mathbb{Z})$ which respects the Hodge
structure.\end{thm}

If $\mathrm{NS}(X)\cong\mathbb{Z}h$ with $h^2=2d=2p_1^{e_1}\ldots
p_m^{e_m}$, where $k\geq 0$, $e_i\geq 1$ and $p_i$ odd primes with
$p_i\neq p_j$ if $i\neq j$, then we consider the Mukai vectors
$$
v^I_J=v^{j_1,\mbox{\ldots},j_s}_{j_{s+1},\mbox{\ldots},j_m}=(p_{j_1}^{e_{j_1}}\ldots
p_{j_s}^{e_{j_s}},h,p_{j_{s+1}}^{e_{j_{s+1}}}\cdots
p_{j_m}^{e_{j_m}}),
$$
where $I=\{j_1, \mbox{\ldots}, j_s\}$ and
$J=\{j_{s+1},\mbox{\ldots} ,j_m\}$ are a partition of
$\{1,\mbox{\ldots},m\}$ such that $I\amalg
J=\{1,\mbox{\ldots},m\}$. The following theorem shows how to
determine $|FM(X)|$ of them corresponding to non-isomorphic moduli
spaces of stable sheaves.

\begin{thm}\label{thm:main1} Let $X$ be a K3 surface with $\mathrm{NS}(X)=\mathbb{Z}h$
such that $h^2=2d=2p_1^{e_1}\ldots p_m^{e_m}$. Then, for all
$v^I_J$ as above, $M_h(v^I_J)$ is a 2-dimensional compact fine
moduli space of stable sheaves on $X$. Moreover, if
$M_h(v^{I_1}_{J_1})\cong M_h(v^{I_2}_{J_2})$, then
$v^{I_1}_{J_1}=v^{I_2}_{J_2}$ or $v^{I_2}_{J_2}=(s_1,h,r_1)$, with
$v^{I_1}_{J_1}=(r_1,h,s_1)$, where the multindexes $I_k$ and $J_k$
vary over all the partitions of $\{1,\mbox{\ldots},m\}$.\end{thm}

\begin{proof} The vectors $v_J^I$ are all isotropic and they are primitive
in $\tilde{H}(X,\mathbb{Z})$, so, by Theorem 5.4 in \cite{Mu2}
$M_h(v^I_J)$ is non-empty. Moreover the hypothesis of Theorem 4.1
in \cite{Mu2}, are satisfied and so the moduli spaces are compact.
By Corollary 0.2 in \cite{Mu3} they are 2-dimensional, while they
are fine by the results in the appendix of \cite{Mu2}.

If $m=0$ or $m=1$ then, by Corollary \ref{cor:num}, we have only
one moduli space with respectively $v=(1,h,1)$ in the first case
and $v=(1,h,p^e)$ in the second case.

Otherwise we must prove that if
$$
v_1=(r_1,h,s_1)=v^{I_1}_{J_1}\neq v^{I_2}_{J_2}=(r_2,h,s_2)=v_2,
$$
with $v_2\neq (s_1,h,r_1)$, then
$$
M_h(v_1)\not \cong M_h(v_2)
$$
But by Theorem \ref{thm:spmod} and Torelli theorem, if we put
$$
M_1:=v_1^\perp/\mathbb{Z}v_1\;\;\;\mbox{and}\;\;\;
M_2:=v_2^\perp/\mathbb{Z}v_2
$$
then it suffices to show that there are no Hodge isometries
between $M_1$ and $M_2$. Obviously, it suffices to show that there
are no Hodge isometries between the transcendental lattices which
lifts to an isometry of the second cohomology groups.

By definition, a representative of a class in $M_i$ ($i=1,2$) is a
vector $(a,b,c)$ such that $bh=as_i+cr_i$, hence
$$bh\equiv as_i
\pmod{r_i},
$$
for $i=1,2$. From now on we will write $(a,b,c)$ for the
equivalence class or for a representative of the class. In fact,
all the arguments we are going to propose are independent from
the choice of a representative.

The Hodge structures on $M_1$ and $M_2$ are induced by the ones
defined on $\tilde{H}(X,\mathbb{Z})$, so, up to an isometry, we
identify $\mathrm{NS}(M_h(v_1))$ and $\mathrm{NS}(M_h(v_2))$ with
$$
S_1:=\langle(0,h,2s_1)\rangle\subset M_1\;\;\;\mbox{and}\;\;\;
S_2:=\langle(0,h,2s_2)\rangle\subset M_2
$$
respectively.

Now, we can describe the transcendental lattices $T_1:=S_1^\perp$
and $T_2:=S_2^\perp$ of $M_h(v_1)$ and $M_h(v_2)$ respectively.

If $(a,b,c)\cdot(0,h,2s_1)=0$ then $bh\equiv 0 \pmod{r_1}$.
Indeed, let us suppose that $bh\equiv K \pmod{r_1}$ where $K\not
\equiv 0 \pmod{r_1}$. Then, by simple calculations, we obtain
$$
(a,b,c)=(L,0,H)+\left(0,n,\frac{n\cdot h-K}{r_1}\right)
$$
as equivalence classes. Here $n=b-kh$, for a particular
$k\in\mathbb{Z}$, $bh\equiv Ls_1 \pmod{r_1}$ and $H$ is an
integer. But now
$$
0=(a,b,c)\cdot (0,h,2s_1)=-2Ls_1+nh=-2Ls_1+(b-kh)h=
$$
\vspace{-0,4cm}
$$
=-2bh+2wr_1+bh-kh^2,
$$
with $w\in\mathbb{Z}$. So
$$
bh\equiv 0 \pmod{r_1}.
$$
This is a contradiction. By these remarks and simple calculations,
a class $y$ in $T_1$, as an element of the quotient $M_1$, has
representative $(0,n,nh/r_1)$. But $(0,n,nh/r_1)\cdot
(0,h,2s_1)=nh=0$. So $y=(0,n,0)$ and
$$
T_1=\{(0,n,0):n\in T_X\}.
$$
Analogously we have
$$
T_2=\{(0,n,0):n\in T_X\}.
$$
By Lemma 4.1 in \cite{Og} (see also point (ii) of Corollary
\ref{cor:num}), if $f:(T_1,\mathbb{C}\omega_1)\rightarrow
(T_2,\mathbb{C}\omega_2)$ is a Hodge isometry, then all the Hodge
isometries from $T_1$ into $T_2$ are $f$ and $-f$. But in this
case $M_1$ and $M_2$ inherit their Hodge structure from
$\tilde{H}(X,\mathbb{Z})$. Hence the two Hodge isometries
$f,g:T_1\rightarrow T_2$ are
$$
(0,n,0) \stackrel{f}{\longmapsto}(0,n,0)\; \mbox{ or }\; (0,n,0)
\stackrel{g}{\longmapsto}(0,-n,0).
$$

Let us show that $f$ cannot be lifted to an isometry from $M_1$
into $M_2$. Equivalently, this means that there are no
isomorphisms between $M_h(v_1)$ and $M_h(v_2)$ which induces $f$.

We start by observing that, if $(a,b,c)\in M_i$ with $i=1,2$, then
$$
(a,b,c)\cdot (0,h,2s_i)\equiv -bh \pmod{r_i}.
$$
Indeed, if $bh\equiv K \pmod{r_i}$ then $a\equiv L \pmod{r_i}$
and so $(a,b,c)=(L,0,H)+(0,n,\frac{n\cdot h-K}{r_i})$ with $n$
and $H$ as before. So, $(a,b,c)\cdot
(0,h,2s_i)=-2bh+2wr_i+bh-kh^2\equiv -bh \pmod{r_i}$ and $nh\equiv
bh \pmod{r_i}$.

Now let us suppose that there is an isometry $\varphi
:M_1\rightarrow M_2$ which induces $f$. We can prove that there is
$(a,b,c)\in M_1$, with $bh\equiv 0 \pmod{r_1}$, such that $\varphi
(a,b,c)=(d,e,f)\in M_2$ with $eh\not \equiv 0 \pmod{r_2}$. First
of all, by our hypotheses about $r_1$ and $r_2$, we can suppose
that there is a prime $p$ which divides $r_2$ but which does not
divide $r_1$ (otherwise we can change the roles of $M_1$ and $M_2$
in the following argument). By Theorem 1.14.4 in \cite{Ni}, there
is an isometry
$$
\psi: H^2(X,\mathbb{Z})\longrightarrow U^3\oplus E_8(-1)^2=\Lambda
$$
such that $k_1:=\psi(h)=(1,d,0,\ldots,0)$, where $h^2=2d$. Let
$k_2:=(0,r_1,0,\ldots,0)$. Now $k_1\cdot k_2=r_1$ and we can take
$n:=\psi^{-1}(k_2)$. Obviously, the vector $(0,n,n\cdot h/r_1)\in
M_1$ is such that $n\cdot h\equiv 0 \pmod{r_1}$. Let us suppose
that $\varphi((0,n,n\cdot h/r_1))=(d,e,f)$ with $e\cdot h\equiv
0\pmod{r_2}$. By the previous remark, this is equivalent to say
that
$$
\varphi\left(\left(0,n,\frac{n\cdot
h}{r_1}\right)\right)=\left(0,m,\frac{m\cdot h}{r_2}\right),
$$
for a given $m\in H^2(X,\mathbb{Z})$.

Since rk$S_1=$rk$S_2=1$, either $\varphi((0,h,2s_1))=(0,h,2s_2)$
or $\varphi((0,h,2s_1))=-(0,h,2s_2)$. In particular, if $\varphi$
correspond to case (1) (the same argument holds if $\varphi$ is as
in case (2)), then
$$
n\cdot h=\left(0,n,\frac{n\cdot
h}{r_1}\right)\cdot\left(0,h,2s_1\right)=\varphi\left(\left(0,n,\frac{n\cdot
h}{r_1}\right)\cdot\left(0,h,2s_1\right)\right)=
$$
\vspace{-0,2cm}
$$
=\left(0,m,\frac{m\cdot
h}{r_2}\right)\cdot\left(0,h,2s_2\right)=m\cdot h.
$$
In particular, $m\cdot h=n\cdot h=r_1$ which is not divisible by
$r_2$. This gives a contradiction and thus $eh\not \equiv 0
\pmod{r_2}$.

The previous remarks show that if
$$
(a,b,c)=\left(0,n,\frac{nh}{r_1}\right)
$$
then
$$
\varphi((a,b,c))=(d,e,f)=(L,0,H)+\left(0,m,\frac{m\cdot
h-K}{r_2}\right),
$$
with $L\not \equiv 0 \pmod{r_2}$. Let us take $(0,N,0)\in T_1$ and
$$
\varphi (0,N,0)=f(0,N,0)=(0,N,0)\in T_2.
$$
Then
$$
(*)\;\;\; nN=\left(0,n,\frac{nh}{r_1}\right)\cdot (0,N,0)=
$$
\vspace{-0,2cm}
$$
=\left[(L,0,H)+\left(0,m,\frac{m\cdot
h-K}{r_i}\right)\right]\cdot (0,N,0)=mN.
$$
Because $H^2(X,\mathbb{Z})$ is unimodular and (*) is true for
every $N\in T_X$, we have $m-n=kh\in \mathrm{NS}(X)$, where
$k\in\mathbb{Z}$. But now $nh=(0,n,\frac{nh}{r_1})\cdot
(0,h,2s_1)=[(L,0,H)+(0,m,\frac{m\cdot h-K}{r_i})]\cdot
(0,h,2s_2)=mh -2Ls_2=nh+kh^2-2Ls_2=nh+2kr_2 s_2-2Ls_2$. So
$L\equiv 0 \pmod{r_2}$ which is contradictory.

Repeating the same arguments for $g$, we see that neither $f$ nor
$g$ lifts to an isometry of the second cohomology groups. So, by
Torelli Theorem, $M_h(v_1)\not \cong M_h(v_2)$.\end{proof}




\section{Genus and polarizations when $\rho=2$}\label{sec:rk2}


In this paragraph we are interested in the number of non
isomorphic FM-partners of K3 surfaces with a given polarization
and Picard number 2.

Our main result is Theorem \ref{thm:main2}. First of all, we
recall the following lemma which is an easy corollary of Nikulin's
Theorem 1.14.2 in \cite{Ni} and whose hypotheses are trivially
verified if $\rho=2$.

\begin{lem}\label{lem:matr} Let $L$ be an even unimodular lattice and let $T_1$ and
$T_2$ be two even sublattice with the same signature $(t_{(+)},
t_{(-)})$, where $t_{(+)}>0$ and $t_{(-)}>0$. Let the
corresponding discriminant groups $(A_{T_1}, q_{T_1})$ and
$(A_{T_2}, q_{T_2})$ be isometric and let $\mathrm{rk}T_1\geq
2+\ell(A_{T_1})$, where $\ell(A_{T_1})$ is the minimal number of
generators of $A_{T_1}$. Then $T_1\cong T_2$.\end{lem}

We prove the following lemma.

\begin{lem}\label{lem:calc} Let $L_{d,n}$ be the lattice $(\mathbb{Z}^2,M_{d,n})$, where
$$
M_{d,n}:=\left(
\begin{array}{cc}
2d & n \\
n & 0
\end{array}
\right),
$$
with $d$ and $n$ positive integers such that $(2d,n)=1$. Then
\\{\rm (i)} the discriminant group $A_ {L_{d,n}}$ is cyclic;
\\{\rm (ii)} if $d_1$, $d_2$,
$n_1$ and $n_2$ are positive integers such that
$(2d_1,n_1)=(2d_2,n_2)=1$ then $A_ {L_{d_1,n_1}}\cong A_
{L_{d_2,n_2}}$ if and only if

{\rm (a.1)} $n_1=n_2$;

{\rm (b.1)} there is an integer $\alpha$ such that $(\alpha,n)=1$
and $d_1\alpha^2\equiv d_2 \pmod{n^2}$;
\\{\rm (iii)} if $L_{d_1,n}\cong L_{d_2,n}$ then one of the following
conditions holds

{\rm (a.2)} $d_1\equiv d_2 \pmod{n}$;

{\rm (b.2)} $d_1d_2\equiv 1 \pmod{n}$.\end{lem}

\begin{proof} Let $e_{d,n}=(1,0)^t$ and $f_{d,n}=(0,1)^t$ be generators of
the lattice $L_{d,n}$. Under the hypothesis $(2d,n)=1$, (i)
follows immediately because
$$
A_{L_{d,n}}:=L_{d,n}^\vee/L_{d,n}
$$
has order $|$det$M_{d,n}|=n^2$ and it is cyclic with generator
$$
\overline{f}_{d,n}:=\frac{ne_{d,n}-2df_{d,n}}{n^2}.
$$
Indeed,
$$
L_{d,n}^\vee=\left\langle\frac{ne_{d,n}-2df_{d,n}}{n^2},\frac{f_{d,n}}{n}\right\rangle
$$
and $\overline{f}_{d,n}$ has order $n^2$ in $A_{L_{d,n}}$.

First we prove that the conditions (a.1) and (b.1) are necessary.
The orders of $A_{L_{d_1,n_1}}$ and $A_{L_{d_2,n_2}}$ are $n_1^2$
and $n_2^2$ respectively with $n_1,n_2>0$, so $n:=n_1=n_2$ (which
is (a.1)). If $A_{L_{d_1,n}}$ and $A_{L_{d_2,n}}$ are isomorphic
as groups, there is an integer $\alpha$ prime with $n$ such that
the isomorphism is determined by
$$
\overline{f}_{d_1,n}\mapsto \alpha\overline{f}_{d_2,n}.
$$
But now
$$
\overline{f}_{d_1,n}^2=\frac{-2d_1}{n^2}\;\;\;\;\;\;\;\;\;\;
\overline{f}_{d_2,n}^2=\frac{-2d_2}{n^2},
$$
and if we want $A_{L_{d_1,n}}$ and $A_{L_{d_2,n}}$ to be isometric
as lattices, we must require
$$
\frac{-2d_1}{n^2}\equiv\alpha^2\frac{-2d_2}{n^2}\pmod{2}.
$$
This is true if and only if
$$
d_1\equiv \alpha^2 d_2 \pmod{n^2}.
$$
So the necessity of condition (b.1) is proved. In the same way it
follows that (a.1) and (b.1) are also sufficient.

Let us consider point (iii). The lattices $L_{d_1,n}$ and
$L_{d_2,n}$ are isometric if and only if there is a matrix
$A\in\mathrm{GL}(2,\mathbb{Z})$ such that
$$
(*)\;\;\;A^tM_{d_1,n}A=M_{d_2,n}.
$$
Let $L_{d_1,n}$ and $L_{d_2,n}$ be isometric and let
$$ A:=\left(
\begin{array}{cc}
x & y \\
z & t
\end{array}
\right).
$$
Then from (*) we obtain the two relations

(1) $d_2=x^2d_1+xzn$;

(2) $2y(yd_1+tn)=0$.

By (2) we have only two possibilities: either $y=0$ or $yd_1=-tn$.
Let $y=0$. From the relation
$$
1=|\mbox{det}(A)|=|xt-yz|=|xt|
$$
it follows that $x=\pm 1$ and so, from (1), we have $d_1\equiv
d_2 \pmod{n}$, which is condition (a.2).

Let us consider the case $yd_1=-tn$. We know that $(d_1,n)=1$ and
hence $y=cn$ and $t=-cd_1$, with $c\in\mathbb{Z}$. From
$$
1=|\mathrm{det}(A)|=|-cxd_1-cnz|
$$
it follows that $c=\pm 1$. We suppose $c=1$ (if $c=-1$ then the
same arguments work by simple changes of signs). Multiplying both
members of relation (1) by $d_1$ we have
$$
d_2d_1\equiv x^2d_1^2 \pmod{n}.
$$
But we know that $\pm 1=\mathrm{det}(A)=-xd_1-nz$ and so
$-xd_1\equiv \pm1 \pmod{n}$. Thus
$$
1\equiv x^2d_1^2 \pmod{n}
$$
and from this we obtain (b.2).\end{proof}

Now we can prove the following theorem (note that point (iii) and
(v) are exactly Theorem 1.7 in \cite{Og}).

\begin{thm}\label{thm:main2} Let $N$ and $d$ be positive integers. Then there are $N$ K3
surfaces $X_1,$\ldots$, X_N$ with Picard number $\rho=2$ such that

{\rm (i)} $X_i$ is elliptic, for every $i\in \{1,\ldots ,N\}$;

{\rm (ii)} there is $i\in \{1,\ldots ,N\}$ such that $X_i$ has a
polarization of degree $2d$;

{\rm (iii)} $\mathrm{NS}(X_i)\not \cong \mathrm{NS}(X_j)$ if $i\neq
j$, where $i,j\in\{1,$\ldots$,N\}$;

{\rm (iv)} $|$det$\mathrm{NS}(X_i)|$ is a square, for every $i\in
\{1,\ldots ,N\}$;

{\rm (v)} there is an Hodge isometry between
$(T_{X_i},\mathbb{C}\omega_{X_i})$ and
$(T_{X_j},\mathbb{C}\omega_{X_j})$, for all
$i,j\in\{1,$\ldots$,N\}$.

\noindent In particular, $X_i$ and $X_j$ are non-isomorphic
FM-partners, for all $i,j\in\{1,$\ldots$,N\}$.\end{thm}

\begin{proof} The surjectivity of the period map for K3 surfaces implies
that, given a sublattice $S$ of $\Lambda$ with rank 2 and
signature $(1,1)$, there is at least one K3 surface $X$ such that
its transcendental lattice $T_X$ is isometric to $T:=S^\perp$.

So the theorem follows if we can show that, for an arbitrary
integer $N$, there are at least $N$ sublattices of $\Lambda$ with
rank 2, signature $(1,1)$  and representing zero which are
non-isometric but whose orthogonal lattices are isometric in
$\Lambda$.

Let us consider in $U\oplus U\hookrightarrow \Lambda$ the
following sublattices
$$
S_{d,n}:=\langle\left(\begin{array}{c} 1 \\ q \\ 0 \\ 0
\end{array}\right),\left(\begin{array}{c} 0 \\ n \\ 1 \\ 0
\end{array}\right)\rangle,
$$
with $(2q,n)=1$ and $n>0$.

We can observe that, when $n$ and $d$ vary, the lattices $S_{q,n}$
are primitive in $\Lambda$ and the matrices associated to their
quadratic forms are exactly the $M_{q,n}$. Since $M_{q,n}$ has
negative determinant, the lattice has signature (1,1). Moreover,
$S_{q,n}$ represents zero

Let $n>2$ be a prime number such that $n>d^2N^4$. We choose
$d_1:=d$, $d_2:=d2^2$,\ldots, $d_{N}=dN^2$. By definition, there
is an integer $\alpha_i$ such that $(\alpha_i,n)=1$ and
$$
\alpha_i^2d_1\equiv d_i \pmod{n^2},
$$
for every $i\in \{1,\ldots,N\}$. Thus the hypotheses (b.1) of
Lemma \ref{lem:calc} are satisfied and by point (ii) of the same
lemma,
$$
A_{S_{d_1,n}}\cong A_{S_{d_i,n}}\cong A_{S_{d_j,n}},
$$
where $i,j\in\{1,...,N\}$. By Lemma \ref{lem:calc}, there are
isometries
$$
\psi_i: S_{d_1,n}^\perp \rightarrow S_{d_i,n}^\perp,
$$
with $i\in\{2,...,N\}$. Now let $(X_1,\varphi_1)$ be a marked K3
surface associated to the lattice $S_{d_1,n}$. By the surjectivity
of the period map we can consider the marked K3 surfaces
$(X_i,\varphi_i)$, with $i\in\{2,...,N\}$, such that

(1)
$\varphi_{i,\mathbb{C}}(\mathbb{C}\omega_{X_i})=\psi_{i,\mathbb{C}}(\varphi_{1,\mathbb{C}}(\mathbb{C}\omega_{X_1}))$;

(2) $\varphi_i(\mathrm{NS}(X_i))=S_{d_i,n}$;

(3) $\varphi_i(T_{X_i})=S_{d_i,n}^\perp$.
\\Obviuosly, the surfaces $X_i$ are FM-partners of $X_1$.

Now we show that, when $i\neq j$,
$$
S_{d_i,n}\not \cong S_{d_j,n}.
$$
First of all we know that, obviously, $d_j\not \equiv d_i
\pmod{n}$ if $i\neq j$. On the other hand,
$$
d_id_j<d^2N^4<n,
$$
so
$$
1\not \equiv d_id_j \pmod{n}.
$$
Hence, by point (iii) of Lemma \ref{lem:calc}, the lattices can
not be isometric. The K3 surfaces $X_1,\ldots, X_N$ are obviously
elliptic and the discriminant of their N\'eron-Severi group is a
square. Moreover $X_1$ has a polarization of degree $2d$.

This shows that it is possible to find $N$ K3 surfaces which
satisfy the hypotheses of the theorem.\end{proof}

The previous theorem gives a new proof of the following result due
to Oguiso (\cite{Og}).

\begin{cor}\label{cor:oguiso} {\bf [13, Theorem 1.7].} Let $N$ be a natural number. Then
there are $N$ K3 surfaces $X_1,$\ldots$, X_N$ with Picard number
$\rho=2$ such that

{\rm (i)} $\mathrm{NS}(X_i)\not \cong \mathrm{NS}(X_j)$ if $i\neq
j$, where $i,j\in\{1,$\ldots$,N\}$;

{\rm (ii)} there is an Hodge isometry between
$(T_{X_i},\mathbb{C}\omega_{X_i})$ and
$(T_{X_j},\mathbb{C}\omega_{X_j})$, for all
$i,j\in\{1,$\ldots$,N\}$.\end{cor}

\begin{remark} The proof proposed by Oguiso in \cite{Og} is based on deep
results in number theory. In particular, it uses a result of
Iwaniec \cite{Iw} about the existence of infinitely many integers
of type $4n^2+1$ which are product of two not necessarily distinct
primes. Theorem \ref{thm:main2} gives an elementary proof of
Theorem 1.7 in \cite{Og} entirely based on simple remarks about
lattices and quadratic forms.\end{remark}

Lemma \ref{lem:matr} is true also when $L=U\oplus U\oplus U$. The
period map is onto also for abelian surfaces (see \cite{Sh}).
Thus, using the lattices $S_{d,n}$ described before, the following
proposition (similar to a result given in \cite{HLOY4}) can be
proved with the same techniques.

\begin{prop}\label{prop:varab} Let $N$ and $d$ be positive integers. Then there are $N$
abelian surfaces $X_1,$\ldots$, X_N$ with Picard number $\rho=2$
such that

{\rm (i)} $\mathrm{NS}(X_i)\not \cong \mathrm{NS}(X_j)$ if $i\neq
j$, with $i,j\in\{1,$\ldots$,N\}$;

{\rm (ii)} there is $i\in\{1,\ldots,N\}$ such that $X_i$ has a
polarization of degree $2d$;

{\rm (iii)} there is an Hodge isometry between
$(T_{X_i},\mathbb{C}\omega_{X_i})$ and
$(T_{X_j},\mathbb{C}\omega_{X_j})$, for all
$i,j\in\{1,$\ldots$,N\}$.\end{prop}

The following easy remark shows that it is possible to obtain an
arbitrarily large number of $M$-polarizations on a $K3$ surface,
for certain $M$.

\begin{remark} Let $N$ be a natural number. Then there are a primitive
sublattice $M$ of $\Lambda$ with signature $(1,0)$ and a $K3$
surface $X$ with $\rho(X)=2$ such that $X$ has at least $N$
non-isomorphic $M$-polarizations. In particular $X$ has at least
$N$ non-isomorphic $M$-polarized FM-partners.

In fact, let $S\cong U$, where $U$ is, as usual, the hyperbolic
lattice. Then, by the surjectivity of the period map, there is a
$K3$ surface $X$ such that $\mathrm{NS}(X)\cong S$.

Let $d$ be a natural number with $d=p_1^{e_1}\ldots p_n^{e_n}$. In
$S$ there are $2^{p(d)-1}$ primitive vectors with autointersection
$2d$. Indeed they are all the vectors of type
$$
f^I_J:=(p_{j_1}^{e_{j_1}}\ldots
p_{j_s}^{e_{j_s}},p_{j_{s+1}}^{e_{j_{s+1}}}\cdots
p_{j_n}^{e_{j_n}}),
$$
for $I$ and $J$ that vary in all possible partitions $I\amalg
J=\{1,\ldots,n\}$.

The group $\mathrm{O}(U)$ has only four elements (i.e. $\pm id$,
the exchange of the vectors of the base and the composition of
this map with $-id$). So it is easy to verify that all these
polarizations are not isomorphic. Choosing $d$ to be divisible by
a sufficiently large number of distinct primes, we can find at
least $N$ non-isomorphic $M$-polarizations. The last assertion
follows from Lemma \ref{lem:matr}.\end{remark}

\medskip

{\small\noindent{\bf Acknowledgements.} The author would like to
express his thanks to Professor Bert van Geemen for his
suggestions and helpful discussions.}

\end{document}